\input ./style/arxiv-general.cfg
\documentclass[MSNbibl,number,citesort,secthm,dvips]{arxbj}
\makeatletter
   \@ifpackageloaded{graphicx}{}{\usepackage{graphicx}}
\makeatother
\usepackage{mathbh}

\volume{22}
\issue{2}
\pubyear{2016}
\firstpage{1278}
\lastpage{1299}
\doi{10.3150/14-BEJ693}
\docsubty{FLA}

\makeatletter
\newcommand{\ds}{\displaystyle}
\newcommand{\eqref}[1]{(\ref{#1})}
\newtheorem{lem}[thm]{Lemma}
\newtheorem{prop}[thm]{Proposition}
\newproclaim{defn}{Definition}[section]
\newproclaim{algo}{Algorithm}[section]
\newremark{rem}{Remark}[section]
\makeatother

\begin{document}
\begin{frontmatter}

\title{Exchangeable exogenous shock models}
\runtitle{Exchangeable exogenous shock models}

\begin{aug}
\author[A]{\inits{J.-F.}\fnms{Jan-Frederik}~\snm{Mai}\thanksref{e1}\ead[label=e1,mark]{jan-frederik.mai@xaia.com}},
\author[A]{\inits{S.}\fnms{Steffen}~\snm{Schenk}\corref{}\thanksref{e2}\ead[label=e2,mark]{steffen.schenk@tum.de}\ead[label=u1,url]{https://www.mathfinance.ma.tum.de/}}
\and
\author[A]{\inits{M.}\fnms{Matthias}~\snm{Scherer}\thanksref{e3}\ead[label=e3,mark]{scherer@tum.de}}
\address[A]{Technische Universit\"at M\"unchen, Parkring 11, 85748
Garching-Hochbr\"uck, Germany.\\ \printead{e1,e2,e3};\\ \printead{u1}}
\end{aug}

%
\received{\smonth{4} \syear{2014}}
%
\revised{\smonth{9} \syear{2014}}

\begin{abstract}
We characterize a comprehensive family of $d$-variate exogenous shock
models. Analytically, we consider a family of multivariate distribution
functions that arises from ordering, idiosyncratically distorting, and
finally multiplying the arguments. Necessary and sufficient conditions
on the involved distortions to yield a multivariate distribution
function are given. Probabilistically, the attainable set of
distribution functions corresponds to a large class of exchangeable
exogenous shock models. Besides, the vector of exceedance times of an
increasing additive stochastic process across independent exponential
trigger variables is shown to constitute an interesting subclass of the
considered distributions and yields a second probabilistic model. The
alternative construction is illustrated in terms of two examples.
\end{abstract}

\end{frontmatter}

\section{Introduction}
\label{secintroductin}
Fatal shock models are standard tools in reliability theory, insurance,
credit risk, and various other fields of application. The present work
characterizes a large class of such models, described by functions
$C\dvtx [0,1]^d\rightarrow[0,1], d\geq2$, of the form
%
\begin{equation}
\label{G}
C(u_1,\ldots,u_d)=\prod
_{k=1}^{d}g_k(u_{(k)}),
\end{equation}
where $g_1=\mbox{id}_{[0,1]}$ is the identity on $[0,1]$, the
mappings $g_k  \dvtx [0,1]\rightarrow[0,1]$ satisfy $g_k(1)=1, k=2,\ldots,d$,
and where $u_{(1)}\leq u_{(2)}\leq\cdots\leq u_{(d)}$ denotes the
ordered list of $u_1,\ldots,u_d$. We derive necessary and sufficient
conditions on $g_k, k=2,\ldots,d$, such that the function $C$ defines a
distribution function (a so-called copula\footnote{A $d$-dimensional
function $C\dvtx [0,1]^d\rightarrow[0,1]$ is a copula if (a) $C$ is the
distribution function of a vector $(U_1,\ldots,U_d)$ on $(\Omega
,\mathcal{F},\mathbb{P})$ and (b) $U_k, k=1,\ldots,d$, is uniformly
distributed on $[0,1]$. For the function $C$ in equation \eqref{G}, (a)
implies (b) as for $(U_1,\ldots,U_d)\sim C$, one has $\mathbb
{P}(U_k\leq u)=C(u,1,\ldots,1)=u$.}), and we construct the
corresponding random vector $\mathbf{X}=(X_1,\ldots,X_d)$ on a
probability space $(\Omega,\mathcal{F},\mathbb{P})$ having $C$ as
distribution function.

It turns out that there is a one-to-one correspondence between copulas
of the form \eqref{G} and random vectors $(X_1,\ldots,X_d)$ with representation
%
\begin{equation}\label{startingpoint}
X_k=\max\bigl\{Z^E\dvt  k\in E\bigr\},\qquad  k=1, \ldots,d, 
\end{equation}
where $Z^E$, $\varnothing\neq E\subseteq\{1,\ldots,d\}$, denote $2^{d-1}$
independent real-valued random variables whose distribution functions
depend only on the cardinality of $E$. More precisely, we show that (a)
any random vector $\mathbf{X}$ defined by \eqref{startingpoint} with
continuous marginal distributions has a (unique) copula of type \eqref{G} and that (b) any copula of the form \eqref{G} represents the
distribution function of a random vector $\mathbf{X}$ that can be
constructed as in \eqref{startingpoint}.

Put differently, as the copula of a random vector $\mathbf{X}$ is the
survival copula of $-\mathbf{X}$, copulas of type \eqref{G}
characterize precisely the set of survival copulas of random vectors
$(\tilde{X}_1,\ldots,\tilde{X}_d)$ with
%
\begin{equation}\label{exchangeable}
\tilde{X}_k:=\min\bigl\{\tilde{Z}^E\dvt  k\in E\bigr\},\qquad
k=1,\ldots,d,
\end{equation}
where the $\tilde{Z}^E$ are identically distributed for subsets $E$
having the same cardinality. Such random vectors $(\tilde{X}_1,\ldots
,\tilde{X}_d)$ are often referred to as exogenous shock models in the
literature (for references see Section~\ref{sec2}), because the $\tilde
{Z}^E$ can be thought of as the arrival times of shocks affecting one
or several constituents in a system of $d$ components, and $\tilde
{X}_k$ represents the first time component $k$ is hit by a shock.

Moreover, with the considered vector $(X_1,\ldots,X_d)$, respectively,
the vector $(\tilde{X}_1,\ldots,\tilde{X}_d)$, being exchangeable, that
is,
\[
(X_1,\ldots,X_d)\stackrel{\mathrm{d}}
{=}(X_{\sigma(1)},\ldots ,X_{\sigma(d)})\qquad  \mbox{for all permutations
$\sigma$ on $\{ 1,\ldots,d\}$},
\]
copulas of type \eqref{G} coincide with the set of survival copulas
corresponding to exchangeable exogenous shock models. Theorem~\ref
{maintheorem} is the main contribution to be derived. For notational
brevity, we introduce a particular set of distribution functions on
$[0,1]$, denoted $\mathcal{D}$ and defined by\footnote{Note that
throughout this document, whenever talking about increasing functions,
we refer to functions that are non-decreasing.}
%
\[
\mathcal{D}:= \bigl\{F\dvtx [0,1]\rightarrow[0,1]\dvt  F\mbox{ continuous and
increasing}, F(1)=1,
F \mbox{ strictly positive on $(0,1]$}\bigr\}.
\]

\begin{thm}\label{maintheorem}
Let $C\dvtx [0,1]^d\rightarrow[0,1]$ have
analytical form \eqref{G}. The following statements are equivalent:
\begin{longlist}[(iii)]
\item[(i)] $C$ is a copula, that is, a multivariate distribution
function.
\item[(ii)] For all $0< u< v\leq1, k\in\mathbb{N}_{0}, j\in\mathbb
{N}$, with $k+j\leq d$, it holds that
\[
G_{j,k}(u,v):=\sum_{i=0}^{j}\pmatrix{j \cr i} (-1)^{i} \prod_{l=1}^{i}g_{l+k}(u)
\prod_{l=i+1}^{j}g_{l+k}(v)\geq0.
\]
\item[(iii)] For all $k\in\mathbb{N}_{0}, j\in\mathbb{N}$, with
$k+j\leq d$, it holds that $H_{j,k}\in\mathcal{D}$, where
\[
H_{j,k}(u):= %
\cases{ \ds\prod_{i=0}^{j-1}g_{k+1+i}^{(-1)^i {j-1 \choose i}}(u),
&\quad $u\in(0,1]$,\vspace*{3pt}
\cr
\ds \lim_{v\searrow0} \prod
_{i=0}^{j-1}g_{k+1+i}^{(-1)^i {j-1 \choose i}}(v), &\quad
$u=0$.} %
\]
\item[(iv)] For all $m\in\{1,\ldots,d\}$, it holds that
$H_{m,d-m}\in\mathcal{D}$.
\end{longlist}
In this case, $C$ is the distribution function of $\mathbf{X}$ as
defined in \eqref{startingpoint}, where $Z^E\sim H_{m,d-m}$ for all
subsets $E$ with cardinality $\vert E\vert=m$. \label{XE}
\end{thm}

\begin{rem}[(Implications of Theorem~\protect\ref{maintheorem})]\label{interpretation}
1. Theorem~\ref{maintheorem} consists of three crucial achievements:
First of all, in (ii), copulas of type \eqref{G} are
characterized purely analytically and it is shown that the functions
$g_k, k=1,\ldots,d$, have to satisfy certain inequality conditions.
Second, in (iii)~respectively (iv), these conditions are
translated to monotonicity requirements, showing that certain
functionals of the $g_k$ have to yield distribution functions on
$[0,1]$. Last but not least, the functionals are interpreted from a
probabilistic point of view by introducing a stochastic model which
manifests the connection to shock models.

2. It is important to note (see Proposition~\ref{every}) that the
$H_{m,d-m}$ defined in Theorem~\ref{maintheorem}(iii) can be
solved for the $g_k$, yielding
\[
g_k=\prod_{m=1}^{d+1-k}H_{m,d-m}^{{d-k \choose m-1}},\qquad
k=1,\ldots,d.
\]
Consequently, $H_{m,d-m}\in\mathcal{D}$, $m=1,\ldots,d$, can be
arbitrary distribution functions, provided that the normalization
constraint
\[
g_1=\prod_{m=1}^{d}H_{m,d-m}^{{d-1 \choose m-1}}=
\mbox{id}_{[0,1]},
\]
which solely stems from the formulation of the theorem in terms of
copulas rather than general multivariate distribution functions, is
fulfilled. Thus, Theorem~\ref{maintheorem} shows that copulas of type
\eqref{G} not only comprise, but precisely consist of the set of
survival copulas corresponding to exchangeable exogenous shock models.

3. The $g_k, k=1,\ldots,d$ can be interpreted as conditional
distribution functions. More precisely, for $u\in(0,1]$ and $(U_1,\ldots
,U_d)\sim C$,
\begin{eqnarray*}
g_k(u)&=&\frac{\prod_{i=1}^{k}g_i(u)}{\prod_{i=1}^{k-1}g_i(u)}=\frac
{\mathbb{P}(U_1\leq u,\ldots,U_{k}\leq u)}{\mathbb{P}(U_1\leq u,\ldots
,U_{k-1}\leq u)}
\\
&=& \mathbb{P}(U_1\leq u,\ldots,U_{k}\leq u\vert
U_1\leq u,\ldots ,U_{k-1}\leq u).
\end{eqnarray*}
Note that for $k\geq3$, $g_k$ corresponds to the ratio between two
diagonal sections of copulas. As such, it is the ratio between a
$k$-Lipschitz and a $(k-1)$-Lipschitz function.
\end{rem}

The remaining sections are organized as follows. Section~\ref{sec2}
gives a brief overview on exogenous shock models and analyzes bivariate
copulas of the form \eqref{G} in more detail. Moreover, a~sketch of the
proof of Theorem~\ref{maintheorem} is given that is thoroughly carried
out in Appendix~\ref{secmain}. Section~\ref{secalternative} provides
an alternative stochastic model for a subclass of type \eqref{G}-copulas that is based on first-passage time constructions involving
additive processes. Furthermore, two examples of corresponding
parametric copula families are given. Finally, Section~\ref{secconclusion} concludes.

\section{Exogenous shock models and bivariate copulas of type \texorpdfstring{\protect\eqref{G}}{(1)}}\label{sec2}
\subsection{Exogenous shock models}\label{secexogenous}
Exogenous shock models have extensively been analyzed in the
literature. The most prominent example is given by the multivariate
distribution (known as Marshall--Olkin distribution) introduced in \cite
{Marshall1967}. The authors consider a random vector $(\tilde
{X}_1,\ldots,\tilde{X}_d)$ as constructed in \eqref{exchangeable},
where the shocks $\tilde{Z}^E$ are exponentially distributed. \cite
{Sarhan2007} relax this condition in the bivariate case by considering
univariate shocks $\tilde{Z}^{\{1\}}, \tilde{Z}^{\{2\}}$ with
generalized exponential distribution. \cite{Kundu2009} extend this idea
and study the law of the vector $(X_1,X_2)$, where $X_1=\max\{Z_1,Z_3\}
$, $X_2=\max\{Z_2,Z_3\}$, and both $Z_1,Z_2$, and $Z_3$ are independent
random variables with generalized exponential distribution. Proceeding
in a similar way, \cite{Li2011} consider a minimum-type construction as
in \eqref{exchangeable} with the shocks $\tilde{Z}^E$ having arbitrary
distribution functions on $(0,\infty)$, and study the resulting
dependence structure in the bivariate case.
Relying on a method of the eponymous authors, \cite{Shoaee2012}
introduce the Block and Basu bivariate generalized exponential
distribution, which results from decomposing the joint distribution
function in its singular and absolutely continuous part and solely
considering the latter. Another direction is pursued by \cite
{Kundu2014}. The authors start with the construction in \eqref
{exchangeable} for $d=2$ and consider Weibull distributed random
variables $\tilde{Z}^E$. Denoting by $F$ the resulting joint
distribution function of $(\tilde{X}_1,\tilde{X}_2)$, they study the
law of $(\tilde{Y}_1,\tilde{Y}_2)$ given by
\[
\tilde{Y}_1=\min\bigl\{\tilde{X}^{\{1,1\}},\ldots,
\tilde{X}^{\{1,N\}}\bigr\},\qquad  \tilde{Y}_2=\min\bigl\{
\tilde{X}^{\{2,1\}},\ldots,\tilde{X}^{\{2,N\}}\bigr\},
\]
where $N$ is a mixing variable with geometric distribution and $(\tilde
{X}^{\{1,n\}},\tilde{X}^{\{2,n\}})_{n\in\mathbb{N}}$ are i.i.d. random
vectors with distribution function $F$.


The ongoing research interest concerning exogenous shock models is
explained by means of their suitability for various applications. \cite
{Marshall1967} show that in their setup with exponentially distributed
shocks, the joint survival function of $(\tilde{X}_1,\tilde{X}_2)$ can
be linked to the first jump times of independent Poisson processes
affecting either one or both components of the random vector. This
derivation is picked up in \cite{Linkskog2003}, where $(\tilde
{X}_1,\tilde{X}_2)$ is identified with the first occurrence time of
windstorm losses in France, respectively, Germany, that are governed by
west (affecting only France), central (affecting only Germany), and
pan-European (affecting both countries) windstorms. The idea of
modelling insurance events by independent counting processes and
mapping the events to one or several claim types is a popular approach
in multivariate insurance models (see the extensive overview in \cite
{Anastasiadis2012}). Alternatively, applying exogenous shock models in
a credit risk context, one might also think of $\tilde{Z}^E$ as arrival
times of economic catastrophes influencing the default of one or
several assets in a portfolio. This idea is applied, for example, in
\cite{Giesecke2003}. In \cite{Baglioni2013}, the authors rely on the
Marshall--Olkin setup (i.e., exponentially distributed shocks) and
model cross-country dependencies between European obligors in order to
compute the systemic risk of the banking sector in the post-crisis era.

While the universe of exogenous shock models in the literature is wide,
most of them either focus on specific distribution functions of the
shocks $\tilde{Z}^E$ (see \cite{Kundu2014a} and the references
therein), or study the vector $(\tilde{X}_1,\ldots,\tilde{X}_d)$ in the
bivariate case, that is, for $d=2$. In the present article, we
consider an arbitrary dimension $d\geq2$ and assume that the
distribution of $\tilde{Z}^E$ solely depends on the cardinality of $E$,
that is, $\tilde{Z}^E$ and $\tilde{Z}^I$ are identically distributed
for $\vert E\vert=\vert I\vert$. Within this setup, Proposition~\ref
{every} shows that the corresponding set of copulas has form \eqref{G}.
For brevity of the proof, we consider the case of continuous, strictly
increasing marginal distribution functions.

\begin{prop}\label{every} Let $\mathbf{X}:=(X_1,\ldots,X_d)$ be defined
as in \eqref{startingpoint}, i.e. $Z^E\stackrel{d}{=}Z^I$ for $\vert
E\vert=\vert I\vert$. Assume that the distribution functions of the
random variables $Z^E$ are continuous and strictly increasing. The
(unique) copula $C$ of $\mathbf{X}$ has the form \eqref{G}.
\end{prop}

\begin{pf}
Let $Z^E$ have the distribution function $F_m$ for $\vert E\vert=m$,
$\varnothing\neq E\subseteq\{1,\ldots,d\}$. The marginal distribution
functions of $\mathbf{X}=(X_1,\ldots,X_d)$ equal
\[
\hat{g}_1(x):=\mathbb{P}(X_k\leq x)=\prod
_{E\dvt k\in E}\mathbb{P}\bigl(Z^E\leq x\bigr)=\prod
_{m=1}^d F_m^{{d-1 \choose m-1}}(x),\qquad
x\in\mathbb{R},
\]
as there are $d-1$ over $m-1$ subsets $E$ with cardinality $m$ and
distribution function $F_{m}$ that appear in the stochastic
construction of $X_k$. By Sklar's theorem, we have to show that the
joint distribution function of $\hat{g}_1(X_1),\ldots,\hat{g}_1(X_d)$
is a copula $C$ of type \eqref{G}. To recognize that, note that
\[
\hat{g}_1(X_k)=\max\bigl\{\hat{g}_1
\bigl(Z^E\bigr)\dvt  k\in E\bigr\},\qquad k=1,\ldots,d,
\]
where $\hat{Z}^E:=\hat{g}_1(Z^E)$ have distribution functions
$H_{m,d-m}:=F_m\circ\hat{g}_1^{-1}\in\mathcal{D}$. As pointed out in
Remark~\ref{interpretation}(2), one can find functions $g_1,\ldots
,g_d$ such that the distribution function of $\hat{g}_1(X_1),\ldots,\hat
{g}_1(X_d)$ is of the form \eqref{G}. Finally, as the required
normalization constraint
\[
g_1:=\prod_{m=1}^{d}
\bigl(F_{m}\circ\hat{g}_1^{-1}
\bigr)^{{d-1 \choose
m-1}}=\mbox{id}_{[0,1]}
\]
is valid, the claim follows like in the proof of Theorem~\ref
{maintheorem}$\mbox{(iv)}\Rightarrow \mbox{(i)}$ in Appendix~\ref{secmain}.
\end{pf}
The crucial condition for a $d$-variate function to represent a copula
is $d$-increasingness (see \cite{Nelsen2006}, Definition~2.10.1, page 43). However, this property is typically non-trivial to
check, with the complexity of the problem increasing exponentially in
higher dimensions.\footnote{This explains why many studies are
restricted to the case $d=2$.} Theorem~\ref{maintheorem}(ii) indicates that for functions of form \eqref{G}, the $d$-increasingness
conditions can be massively simplified and reduce to the verification
of $G_{j,k}(u,v)\geq0$ for certain indices $k\in\mathbb{N}_0, j\in
\mathbb{N}$, and certain pairs $(u,v)\in[0,1]^2$. To get accustomed to
the paper and to develop a deeper understanding of these conditions,
the following section analyzes copulas of type \eqref{G} in the simpler
bivariate case.

\subsection{Bivariate copulas of type \texorpdfstring{\protect\eqref{G}}{(1)}}
In \cite{Durante2007}, the authors study copulas of type \eqref{G} for
the special case $g_2=g_3=\cdots=g_d$. In the bivariate case, such
copulas coincide with the more general class considered in the present
article. By Theorem~\ref{maintheorem}, $C\dvtx [0,1]^2\rightarrow[0,1]$,
$C(u_1,u_2)=g_1(u_{(1)}) g_2(u_{(2)})$, is a copula if and only if
$G_{j,k}(u,v)\geq0$ for $(k,j)\in\{(0,1), (0,2), (1,1)\}$ and $u,v\in
[0,1], u<v$, that is, if and only if $g_1(v)-g_1(u)\geq0$,
$g_2(v)-g_2(u)\geq0$, and $g_1(v) g_2(v)-2 g_1(u) g_2(v)+g_1(u)
g_2(u)\geq0$.

While the first two conditions are easy to interpret and imply that
$g_1$ and $g_2$ have to be increasing (note that $g_1=\mbox{id}_{[0,1]}$ is increasing by definition already), the third one is
more interesting to analyze. It is shown in \cite{Durante2008} and the
proof of the main theorem below that, given increasingness of $g_1$ and
$g_2$, the third condition is equivalent to $g_2$ being strictly
positive and continuous and $g_1/g_2$ being increasing on $(0,1]$. It
is easy to verify that $C$ is the copula of $\mathbf{X}=(X_1,X_2)$, where
\begin{eqnarray*}
X_1&= & \max\bigl\{Z^{\{1\}},Z^{\{1,2\}}\bigr\},
\\
X_2&=& \max\bigl\{Z^{\{2\}},Z^{\{1,2\}}\bigr\},
\end{eqnarray*}
and $Z^{\{1\}}, Z^{\{2\}}, Z^{\{1,2\}}$ are independent random
variables with distribution functions\footnote{Strictly speaking, in
order for $g_2$ and $g_1/g_2$ to be proper distribution functions on
$[0,1]$, one has to consider their right-continuous extensions at zero
as considered in the definition of $H_{j,k}$ in Theorem~\ref{maintheorem}(iii).} $Z^{\{1\}}, Z^{\{2\}}\sim g_2$ and $Z^{\{
1,2\}}\sim g_1/g_2$.

Several dependence properties of $C$ can be derived in closed form. For
instance, the lower and upper tail dependence coefficients $\lambda_L$
and $\lambda_U$ equal
\begin{eqnarray*}
\lambda_L &:= & \lim_{u\searrow0}\frac{C(u,u)}{u}=\lim
_{u\searrow0} g_2(u),
\\
\lambda_U &:=& \lim_{u\nearrow1}\frac{C(u,u)-2 u+1}{1-u}=1-g_2^{\prime}(1-),
\end{eqnarray*}
where $g_2^{\prime}(1-)$ denotes the left-sided derivative of $g_2$ at
$u=1$, which exists by monotonicity of~$g_2$. Further dependence
properties are derived in \cite{Durante2006,Durante2007}, including
measures of association and extremal dependence coefficients.
As an example for $C$, define $g_2 \dvtx [0,1]\rightarrow[0,1]$ by
\[
g_2(u):=\min\{a u+b,1-c+c u\}, \qquad  a>1, b>0, c>0, b+c\leq1.
\]
Specified in that way, $g_2$ starts at $b$ with slope $a$, has a kink
at $u=(1-(b+c))/(a-c)$, continues to increase with slope $c$, and ends
at $g_2(1)=1$. Applying the tail dependence formulas, it follows that
$\lambda_L=\lim_{u\searrow0} g_2(u)=b$ and $\lambda
_U=1-g_2^{\prime}(1-)=1-c$. Thus, copulas of type \eqref{G} can admit both
positive upper and lower tail dependence with arbitrary values in $[0,1]$.

\subsection{Extreme-value copulas of type \texorpdfstring{\protect\eqref{G}}{(1)}}
One may determine the intersection between copulas having form \eqref
{G} and extreme-value copulas. A $d$-variate extreme-value copula $C$ satisfies
%
\begin{equation}\label{exValue}
C\bigl(u_1^t,\ldots,u_d^t
\bigr)=C^t(u_1,\ldots,u_d)\qquad  \mbox{for
all $t>0, u_1,\ldots,u_d\in[0,1]$}.
\end{equation}
Proposition~\ref{propex} shows that extreme-value copulas of type
\eqref{G} correspond to choosing power functions for $g_k$, where the
sequence of exponents must be $d$-monotone.

\begin{defn}[($d$-monotone sequence)] A real-valued sequence $\{a_0,\ldots
,a_{d-1}\}$ is called $d$-monotone if
\[
\sum_{i=0}^{j-1}(-1)^i \pmatrix{j-1 \cr i}a_{k+i}\geq0 \qquad \mbox {for all $k\in\mathbb{N}_0,
j\in\mathbb{N}\dvt  k+j\leq d$}.
\]
\end{defn}

\begin{prop}[(Extreme-value copulas of type \eqref{G})]\label{propex}
Let $C$ have the form \eqref{G}.
$C$ is an extreme-value copula if and only if $g_k(u)=u^{a_{k-1}},
k=1,\ldots,d$, for $u\in(0,1]$ and a $d$-monotone sequence $\{a_0,\ldots
,a_{d-1}\}$ with $a_0=1$.
\end{prop}

For a proof, see Appendix~\ref{appextremevalue}. The corresponding
class of extreme-value copulas is well-known in the literature: It is
precisely the exchangeable family of Marshall--Olkin survival copulas
(see, e.g., \cite{Ressel2013}). In Section~\ref{secexamples}, we will
investigate this example in more detail.

\subsection{Strategy to prove Theorem \texorpdfstring{\protect\ref{maintheorem}}{1.1}}
Apart from some technical lemmata, the proof of Theorem~\ref
{maintheorem} provides valuable insights into the structure of the
objects $G_{j,k}$ and $H_{j,k}$ and their relation to the stochastic
model in equation \eqref{startingpoint}. We are going to show that
$\mbox{(iv})\Rightarrow\mbox{(i)}\Rightarrow\mbox{(ii)}\Rightarrow\mbox{(iii)}\Rightarrow\mbox{(iv)}$. The
central ideas can be summarized as follows. The rigorous proof of
Theorem~\ref{maintheorem} is given in Appendix~\ref{secmain}.

\begin{rem}[(Structure of the proof of Theorem~\ref{maintheorem})]
$\mbox{(iv)}\Rightarrow \mbox{(i)}$  Starting with the random vector
$\mathbf{X}=(X_1,\ldots,X_d)$ given in the theorem, one can compute
that each $X_k, k=1,\ldots,d$, is uniformly distributed on $[0,1]$ and
that $C$ is the distribution function of $\mathbf{X}$, hence a copula.

$\mbox{(i)}\Rightarrow\mbox{(ii)}$  Being a copula, $C$ induces a
probability measure $dC$ on $[0,1]^d$. It can be deduced that
$G_{j,k}(u,v)$ corresponds to the mass assigned by $dC$ to certain
subsets of $[0,1]^d$. Therefore, it has to be greater than or equal to zero.

$\mbox{(ii)}\Rightarrow\mbox{(iii)}$  This is the most difficult and
lengthy part of the proof. Besides minor technical conditions, the
central task is to show that non-negativity of $G_{j,k}(u,v), u<v$,
implies increasingness of $H_{j,k}$. The underlying proof idea is to
split up $G_{j,k}(u,v)$ into two summands, one involving the difference
$H_{j,k}(v)-H_{j,k}(u)$, the other one corresponding\vspace*{1pt} to the probability
mass $dC(I)$ induced by a copula $C$ of type \eqref{G} for a subset
$I\subset[0,1]^d$. For a sufficiently ``small'' subset $I$, it is shown
that the sign of $G_{j,k}(u,v)$ is dominated by the first part, that
is, the difference $H_{j,k}(v)-H_{j,k}(u)$. Thus, for $G_{j,k}(u,v)$
to be non-negative, $H_{j,k}(v)-H_{j,k}(u)$ has to be greater than or
equal to zero, which establishes the claimed increasingness of $H_{j,k}$.

$\mbox{(iii)}\Rightarrow\mbox{(iv)}$  This is trivial as $\mbox{(iv)}$ is a special case of (iii).
\end{rem}

%

\section{Alternative construction via additive processes}\label{secalternative}
\subsection{Additive frailty construction}\label{secconstruction}
Theorem~\ref{maintheorem} has shown that any $d$-dimensional copula of
the form \eqref{G} arises from the stochastic construction in \eqref
{startingpoint}, involving $2^d-1$ random variables $Z^E, \varnothing\neq
E\subseteq\{1,\ldots,d\}$. Consequently, due to the exponentially
increasing number of random objects that have to be sampled, simulation
becomes practically impossible in large dimensions. This section
provides an alternative construction for a subclass of type \eqref
{G}-copulas based on a first-passage time construction with additive processes.

According to \cite{Sato2007}, Definition~1.6, page~3, a stochastic
process $\{\Lambda_t\}_{t\geq0}$ on $\mathbb{R}$ is called additive if
it starts at zero almost surely, is stochastically continuous, admits
independent increments, and is c\`adl\`ag almost surely. Increasing
additive processes, called additive subordinators in the sequel, are
closely connected with Bernstein functions. A function $\Psi \dvtx [0,\infty
)\rightarrow[0,\infty)$ is called Bernstein function\footnote{For more
details on Bernstein function, see \cite{Schilling2010}.} if $\Psi
(0)=0$, $\Psi$ is infinitely often differentiable on $(0,\infty)$ with
$(-1)^{n-1} \Psi^{(n)}(x)\geq0$ for all $n\in\mathbb{N}$ and $x>0$.
The law of an additive subordinator $\{\Lambda_t\}_{t\geq0}$ can be
described by a family $\{\Psi_t\}_{t\geq0}$ of Bernstein functions
subject to certain consistency conditions (see (i)--(iii) below).
Denoting the Laplace transform operator by $\mathcal{L}$, it follows
from \cite{Sato2007}, page~47 ff., that
\begin{longlist}[(iii)]
\item[(i)] $\Psi_0(x)\equiv0$ for all $x\geq0$,
\item[(ii)] $\Psi_t-\Psi_s$ is a Bernstein function for all $0\leq s\leq t$,
\item[(iii)]  $\mathcal{L} (\Lambda_t-\Lambda_s)=\exp(-\Psi_t+\Psi
_s)$ for all $0\leq s\leq t$.
\end{longlist}
It is well known that $\Psi_t$ admits a L\'evy--Khintchine
representation, that is,
%
\begin{equation}\label{Khinthine}
\Psi_t(x)=a_t \mathbh{1}_{\{x>0\}}+b_t
x+\int_{(0,\infty)}\bigl(1-\mathrm{e}^{-x
s}\bigr)
\nu_t(\mathrm{d}s), \qquad x\geq0,
\end{equation}
with a L\'evy measure $\nu_t$ on $(0,\infty)$ and parameters
$a_t,b_t\geq0$.
One of the most prominent examples for increasing additive processes
are L\'evy subordinators, which not only exhibit independent, but even
stationary increments. For a L\'evy subordinator, the corresponding
family $\{\Psi_t\}_{t\geq0}$ satisfies $\Psi_t=t \Psi_1, t\geq0$.
Another example is given by so-called self-similar additive
subordinators (see \cite{Sato2007}, page 99 ff.), which are
characterized by $\Psi_t(x)=\Psi_1(x t^H), x,t\geq0$, for an $H>0$ and
a specific Bernstein function $\Psi_1$.
%
%

Now consider an additive subordinator $\Lambda=\{\Lambda_t\}_{t\geq0}$
with $\lim_{t\rightarrow\infty}\Lambda_t=\infty$ and define a sequence
$\{X_k\}_{k\in\mathbb{N}}$ of random variables by
%
\begin{equation}\label{firstpassage}
X_k:=\inf\{t\geq0\dvt  \Lambda_t\geq E_k\}, \qquad k
\in\mathbb{N},
\end{equation}
where $E_k, k\in\mathbb{N}$, are i.i.d. unit exponentially distributed
random variables that are independent of $\{\Lambda_t\}_{t \geq0}$. By
construction, $\{X_k\}_{k\in\mathbb{N}}$ is an exchangeable sequence of
random variables. The following proposition outlines that
the (unique) survival copula of $(X_1,\ldots,X_d)$, denoted $C_{\Lambda
,d}$ in the sequel, is of type \eqref{G} for any $d\geq2$.

\begin{prop}\label{subclass}
Define a sequence $\{X_k\}_{k\in\mathbb
{N}}$ of random variables as in \eqref{firstpassage}. Let $\{\Psi_t\}
_{t\geq0}$ be the family of Bernstein functions corresponding to the
increasing additive process $\Lambda=\{\Lambda_t\}_{t\geq0}$ and denote
by $\bar{F}_1$ the survival function of $X_1$. The survival copula
$C_{\Lambda,d}$ of $(X_1,\ldots,X_d)$ has the form~\eqref{G} for any
$d\geq2$, with
\[
g_k(u):=\exp \bigl(-\Psi_{\bar{F}_1^{-1}(u)}(k)+\Psi_{\bar
{F}_1^{-1}(u)}(k-1)
\bigr), \qquad k=1,\ldots,d.
\]
\end{prop}

The proof is to be found in Appendix~\ref{appsubclass}. Referring to
the numerical motivation for the first-passage time construction at the
beginning of this section, a generic sampling algorithm for $C_{\Lambda
,d}$ can be stated as follows.

\begin{algo}[(Simulation of $C_{\Lambda,d}$ in Proposition~\ref
{subclass})]\label{algo2}
1. Simulate $d$ independent, unit exponentially distributed random
variables $E_1,\ldots,E_d$.

2. Simulate one path of $\{\Lambda_t\}_{t\geq0}$
until $\Lambda_t\geq\max\{E_1,\ldots,E_d\}$.

3. Compute $X_k:=\inf\{t\geq0\dvt \Lambda_t\geq E_k\}, k=1,\ldots,d$.

4. Set $U_k:=\bar{F}_1(X_k), k=1,\ldots,d$, and return $(U_1,\ldots,U_d)$.
\end{algo}

Clearly, the central task is to simulate the path of the additive
process $\Lambda$. Provided this can be accomplished efficiently, the
algorithm provides a fast sampling routine even in large dimensions.
For L\'evy processes of (compound) Poisson type, see \cite{Sato2007}, page 17 ff. The path generation of self-similar additive processes
can be accomplished via more general results in \cite{marsaglia1984fast}.
Last but not least, the construction of Dirichlet processes (which might be viewed as an elementary
transform of a special class of additive subordinators) is discussed in
\cite{Fe74}.

Extending $(X_1,\ldots,X_d)$ to larger dimensions solely requires the
simulation of further i.i.d. exponentially distributed triggers
$E_{d+1}, E_{d+2},\ldots,$ possibly supplemented by simulating
additional increments of $\Lambda$ until the largest trigger is
exceeded (it can be shown that $\lim_{d\rightarrow\infty}\mathbb
{E}[E_{(d)}]/\log d=1$). Besides, we consider it interesting to study
the coherence between properties of the additive process in the
first-passage time setup \eqref{firstpassage} and the distribution
functions of $Z^E$, respectively, $\tilde{Z}^E$, in the maximum
construction in \eqref{startingpoint}, respectively, \eqref{exchangeable}.

\subsection{Examples of tractable families}\label{secexamples}

\begin{longlist}[(2)]
\item[(1)] \textit{Extendible Marshall--Olkin copulas}

The first example assigns a very special meaning to the survival copula
of the random vector $(X_1,\ldots,X_d)$ defined in \eqref{firstpassage}
when choosing $\{\Lambda_t\}_{t\geq0}$ to be a L\'evy subordinator. If
$C$ in \eqref{G} is given by the functions $g_k(u)=u^{a_{k-1}},
k=1,\ldots,d$, and defines a copula for any $d\geq2$, it can be
constructed in two quite different ways. On the one hand, $C$ arises as
the survival copula of the additive frailty construction when plugging
in a L\'evy subordinator for $\{\Lambda_t\}_{t\geq0}$. On the other
hand, $C$ is the survival copula of $(\tilde{X}_1,\ldots,\tilde{X}_d)$
in \eqref{exchangeable} with exponentially distributed shocks $\tilde
{Z}^E$. Interestingly, it can be deduced that the stationary increments
of the L\'evy process in the frailty setup translate to the
characterizing lack-of-memory property of the exponential distribution
in the shock construction.

To recognize this coherence, consider a copula $C$ of type \eqref{G}
with $g_k(u)=u^{a_{k-1}}, k=1,\ldots,d$, for a real-valued sequence $\{a_0,\ldots,a_{d-1}\}$ with $a_0=1$. In \cite{Mai2011}, Theorem~2.3, the
authors show that $C$ corresponds to the exchangeable subclass of the
multivariate distribution function introduced in \cite{Marshall1967},
which is why $C$ is called an exchangeable Marshall--Olkin copula. Put
differently, $C$ is the survival copula of $(\tilde{X}_1,\ldots,\tilde
{X}_d)$ in \eqref{exchangeable}, where the shocks $\tilde{Z}^E$ are
exponentially distributed, with rates depending on the cardinality of $E$.

By Theorem~\ref{maintheorem}, a function $C$ of type \eqref{G} with
$g_k(u)=u^{a_{k-1}}, k=1,\ldots,d$, and a real-valued sequence $\{
a_0,\ldots,a_{d-1}\}$, $a_0=1$, defines a copula for any $d\geq2$ if
and only if $H_{j,k}\in\mathcal{D}$ for all $j\in\mathbb{N}, k\in\mathbb
{N}_0$. As discussed in the proof of Proposition~\ref{propex}, this is
naturally equivalent to $\{a_0,\ldots,a_{d-1}\}$ being $d$-monotone for
any $d\geq2$. A sequence $\{a_k\}_{k\in\mathbb{N}_0}$ that is
$d$-monotone for any $d\geq2$ is called completely monotone. Combining
\cite{Mai2009a}, Lemma~2.6 with \cite{Gnedin2008}, Corollary~4.2, $\{
a_k\}_{k\in\mathbb{N}_0}$ with $a_0=1$ is completely monotone if and
only if there exists a L\'evy subordinator $\Lambda=\{\Lambda_t\}_{t\geq
0}$ characterized by $\{\Psi_t\}_{t\geq0}, \Psi_t=t \Psi_1$, such that
$\Psi_1(1)=1$ and $a_{k-1}=\Psi_1(k)-\Psi_1(k-1)$ for all $k\in\mathbb
{N}$. Applying Proposition~\ref{subclass}, $C_{\Lambda,d}$ is a
survival copula of type \eqref{G} for all $d\geq2$, with
\begin{eqnarray*}
g_k(u)&:=& \exp \bigl(-\Psi_{\bar{F}_1^{-1}(u)}(k)+\Psi_{\bar
{F}_1^{-1}(u)}(k-1)
\bigr), \qquad k\in\mathbb{N},
\\
\bar{F}_1(x)&=& \mathrm{e}^{-\Psi_x(1)}\stackrel{(\Psi_x=x
\Psi_1)}{=}\mathrm{e}^{-x \Psi
_1(1)}\stackrel{\bigl(\Psi_1(1)=1
\bigr)} {=}\mathrm{e}^{-x}.
\end{eqnarray*}
Substituting $\bar{F}_1^{-1}(u)=-\log(u)$ in $g_k$, it follows that
\begin{eqnarray*}
g_k(u)&=&\exp \bigl(-\Psi_{-\log(u)}(k)+\Psi_{-\log(u)}(k-1)
\bigr)
\\
&=&\exp \bigl(\log(u) \Psi_1(k)-\log(u)\Psi_1(k-1)
\bigr)=u^{\Psi_1(k)-\Psi
_1(k-1)}=u^{a_{k-1}}.
\end{eqnarray*}
%

\item[(2)]  \textit{Sato-frailty copulas}

The second example concerns the survival copula of $(X_1,\ldots,X_d)$
in \eqref{firstpassage} for self-similar additive processes $\{\Lambda
_t\}_{t\geq0}$, which admits a nice simple form and might be attractive
for applications and further analyses. A self-similar additive
subordinator $\Lambda=\{\Lambda_t\}_{t\geq0}$, sometimes referred to as
increasing Sato process in the literature, can be identified with a
family $\{\Psi_t\}_{t\geq0}$ of Bernstein functions satisfying $\Psi
_t(x)=\Psi_1(x t^H), x,t\geq0$, for an $H>0$ and a so-called
self-decomposable Bernstein function $\Psi_1$. The attribute
``self-decomposable'' means that $\Psi_1$ admits a representation as in
\eqref{Khinthine} with $a_1=0$ and with $\nu_1$ possessing a density
$\nu_1(\mathrm{d}s)=k(s)/s \,\mathrm{d}s$ such that $s\mapsto k(s)$ is decreasing on
$(0,\infty)$. By Proposition~\ref{subclass}, the copula $C_{\Lambda,d}$
has the form \eqref{G} for all $d\geq2$, with $g_k, k=1,\ldots,d$,
given by
\begin{eqnarray*}
g_k(u)&:=& \exp \bigl(-\Psi_{\bar{F}_1^{-1}(u)}(k)+\Psi_{\bar
{F}_1^{-1}(u)}(k-1)
\bigr),\qquad  k\in\mathbb{N},
\\
\bar{F}_1(x)&=& \mathrm{e}^{-\Psi_x(1)}\stackrel{(\Psi_x(1)=
\Psi _1(x^H))} {=}\mathrm{e}^{-\Psi_1(x^H)}.
\end{eqnarray*}
Plugging $\bar{F}^{-1}(u)= (\Psi_1^{-1}(-\log u) )^{1/H}$ into
$g_k$ yields
\begin{eqnarray*}
g_k(u)&=&\exp \bigl(-\Psi_{ (\Psi^{-1}(-\log u) )^{1/H}}(k)+\Psi
_{ (\Psi^{-1}(-\log u) )^{1/H}}(k-1) \bigr)
\\
&=&\exp \bigl(-\Psi_1 \bigl(k \Psi_1^{-1}(-
\log u) \bigr)+\Psi_1 \bigl((k-1) \Psi_1^{-1}(-
\log u) \bigr) \bigr)
\\
&=&\frac{\psi (k \psi^{-1}(u) )}{\psi ((k-1) \psi^{-1}(u) )}
\end{eqnarray*}
for $\psi\dvtx [0,\infty)\rightarrow(0,1], \psi(x):=\exp(-\Psi_1(x))$. The
resulting class $C_{\Lambda,d}$ is analytically tractable and is
analyzed in detail and illustrated in \cite{Mai2014}. As a concluding
example, consider the self-decomposable Bernstein function
\[
\Psi_1(x):=\int_{0}^{\infty}
\bigl(1-\mathrm{e}^{-x t}\bigr) \underbrace{\beta \frac{\exp
(-\eta t)}{t}
\,\mathrm{d}t}_{:=\nu(\mathrm{d}t)}=\beta \log \biggl(1+\frac{x}{\eta} \biggr),\qquad x,\beta,
\eta>0,
\]
which obviously satisfies the required property that $\nu(\mathrm{d}s)=k(s)/s
\,\mathrm{d}s$ for a decreasing function $s\mapsto k(s)$. Defining $\psi \dvtx  [0,\infty
)\rightarrow(0,1]$, $\psi(x):=\exp(-\Psi_1(x))$, and choosing $\Lambda
=\{\Lambda_t\}_{t\geq0}$ to be an increasing Sato process with $\Psi
_t(x)=\Psi_1(x t^H), x,t\geq0$, for an $H>0$, yields
\begin{eqnarray*}
C_{\Lambda,d}(u_1,\ldots,u_d)&=&\prod
_{k=1}^{d} \frac{\psi (k \psi
^{-1}(u_{(k)}) )}{\psi ((k-1) \psi^{-1}(u_{(k)}) )}
\\[-1.5pt]
&=& \Biggl(\prod_{k=1}^d
\frac{1+(k-1)  (u_{(k)}^{-{1}/{\beta
}}-1 )}{1+k  (u_{(k)}^{-{1}/{\beta}}-1 )} \Biggr)^\beta,
\end{eqnarray*}
which in the bivariate case admits the form $C_{\psi
,2}(u_1,u_2)=u_{(1)}/ (2-u_{(2)}^{1/\beta} )^\beta$.
Summarizing,\vspace*{1pt} the construction in \eqref{firstpassage} can be used to
derive interesting new parametric families of\vspace*{-3pt} copulas.
\end{longlist}

\section{Conclusion}\label{secconclusion}
%
The main object of our study are copulas of functional form \eqref{G},
which are parameterized by $d$ functions $g_k, k=1,\ldots,d$. Three
major findings constitute the core of the present article. First,
necessary and sufficient conditions in terms of inequalities for the
$g_k$ are derived. Second, these conditions are shown to relate to
monotonicity conditions for multiplicative conjunctions of the~$g_k$.
Third, the monotonicity restrictions are embedded into a stochastic
model for the corresponding copula. As a result, it turns out that the
considered class of copulas is interrelated with the set of
exchangeable exogenous shock models, which have been analyzed in the
literature on the level of specific models, however, to the best of our
knowledge not on a unified level as in the present work. In addition,
it is outlined how a subclass of the considered set of copulas can
alternatively be constructed via additive processes. The latter finding
seems to be promising in two regards. On the one hand, it provides an
alternative sampling approach that may be beneficial for
high-dimensional simulation purposes. On the other hand, by combining
the characterization properties in Theorem~\ref{maintheorem} with the
alternative construction, interesting theoretical results -- as
illustrated by the Marshall--Olkin example in Section~\ref{secexamples} -- may be\vspace*{-4.5pt} derived.

\begin{appendix}
\section{}\label{appextremevalue}\vspace*{-3pt}
\begin{pf*}{Proof of Proposition~\protect\ref{propex}}
If $g_k(u)=u^{a_{k-1}}, k=1,\ldots,d$, for a sequence $\{a_0,\ldots
,a_{d-1}\}$, the functions $H_{j,k}$ in Theorem~\ref{maintheorem}(iii) are given by
\[
H_{j,k}(u) =\prod_{i=0}^{j-1}g_{k+1+i}^{(-1)^i {j-1 \choose
i}}(u)=
\prod_{i=0}^{j-1}u^{(-1)^i {j-1 \choose i} a_{k+i}}=
u^{\sum
_{i=0}^{j-1}(-1)^i {j-1 \choose i} a_{k+i}}.
\]
Thus, the $H_{j,k}$ are distribution functions in $\mathcal{D}$ (i.e.
$C$ is a copula) if and only if the sequence $\{a_0,\ldots,a_{d-1}\}$
is $d$-monotone. Moreover, it is apparent that due to the power
function structure of the $g_k$, $C$ satisfies the extreme-value property.

It remains to show that any extreme-value copula of type \eqref{G}
implies a power function structure for the $g_k$. By setting
$u_1=u_2=u\in[0,1]$ and $u_3=\cdots=u_d=1$ in \eqref{exValue},
%
\[
g_2\bigl(u^t\bigr)= g_2(u)^t
\qquad \mbox{for all $u\in(0,1], t>0$.}
\]
Defining $\Theta \dvtx  [0,\infty)
\rightarrow\mathbb{R}, \Theta (x):=g_2(\exp(-x))$, this is
equivalent to
\begin{equation}\label{ex1}
\Theta(t x)=\Theta(x)^t \qquad \mbox{for all $x,t>0$.}
\end{equation}
Thus, setting $t=n\in\mathbb{N}$ and $x=1/n$, it holds that
\begin{equation}\label{ex2}
\Theta(1)=\Theta \biggl(\frac{1}{n} \biggr)^n \quad \Longrightarrow\quad
\Theta \biggl(\frac
{1}{n} \biggr)=\Theta(1)^{{1}/{n}}.
\end{equation}
Consequently, for all $x\in\mathbb{Q}\cap(0,\infty)$, $x=p/q, p,q\in
\mathbb{N}$, we have
\[
\Theta \biggl(\frac{p}{q} \biggr)\stackrel{\mathrm{\protect\eqref{ex1}}} {=}\Theta \biggl(
\frac
{1}{q} \biggr)^p\stackrel{\mathrm{\eqref{ex2}}} {=}
\Theta(1)^{{p}/{q}}.
\]
By Theorem~\ref{maintheorem}(iii), it holds that $H_{1,1}=g_2\in
\mathcal{D}$, implying that $g_2(u)>0$ for $u\in(0,1]$. Therefore, one
deduces that $0<\Theta(1)=g_2(\exp(-1))\leq1$ and the previous equation
yields $\Theta(x)=\exp(-a_{1} x)$ for all $x\in\mathbb{Q}\cap(0,\infty
)$, where $a_1=-\log (\Theta(1) )\geq0$. As $\mathbb{Q}\cap
(0,\infty)$ is dense in $\mathbb{R}^+$, it follows that
\[
g_2 \bigl(\mathrm{e}^{-x} \bigr)=\Theta(x)=\mathrm{e}^{-a_{1} x}=
\bigl(\mathrm{e}^{-x} \bigr)^{a_1},\qquad  x\in\mathbb{R}^+.
\]
By continuity of $g_2$ at zero, $g_2(u)=u^{a_1}$ for all $u\in[0,1]$.
Iteratively, by setting $u_1=\cdots=u_l=u\in[0,1]$, $u_{l+1}=\cdots
=u_d=1$, and subsequently raising $l$, we conclude that
$g_k(u)=u^{a_{k-1}}, k=1,\ldots,d$, for parameters $a_0,\ldots
,a_{d-1}\geq0, a_0=1$.\footnote{An alternative way to constitute the
power function structure of the $g_k$ for extreme-value copulas of type
\eqref{G} is via Pickand's theorem as considered in \cite{Durante2010}.
Comparing the Pickands representation of $C$ with the functional form
in \eqref{G} on the diagonal $u=u_1=\cdots=u_d$, it can be shown by
induction that $g_k(u)=u^{a_{k-1}}$ for a parameter
\[
a_{k-1}=k P \biggl(\underbrace{\frac{1}{k},\ldots,
\frac{1}{k}}_{k \,\mathrm{times}},0,\ldots,0 \biggr)-(k-1) P \biggl(\underbrace{
\frac
{1}{k-1},\ldots,\frac{1}{k-1}}_{(k-1) \ \mathrm{times}},0,\ldots ,0 \biggr)
\]
depending on the Pickand dependence function $P$ for fixed values. By
the first part of the proof, the claim follows. We thank the referee
for pointing us to the idea for this alternative proof.
}
\end{pf*}

\section{}\label{secmain}
\begin{pf*}{Proof of Theorem~\protect\ref{maintheorem}}
$\mbox{(iv)}\Rightarrow\mbox{(i)}$: If $\mbox{(iv)}$ holds, the functions
$H_{m,d-m}$ defined in Theorem~\ref{maintheorem}(iii) are valid
distribution functions on $[0,1]$, and we can consider the
corresponding stochastic model given in the theorem. We will show that
the resulting distribution function of $(X_1,\ldots,X_d)$ is a copula
given by $C$. First of all, recognize that each $X_k, k=1,\ldots,d$,
has a uniform marginal distribution due to
\begin{eqnarray*}
\mathbb{P}(X_k\leq u)&=& \prod_{E\dvt  k\in E}
\mathbb{P} \bigl(Z^E\leq u \bigr)\stackrel{(\ast)} {=}\prod
_{m=1}^{d} H_{m,d-m}^{{d-1 \choose m-1}}(u)
\\
&=& \prod_{m=1}^{d} \Biggl(\prod
_{i=0}^{m-1}g_{d-m+1+i}^{(-1)^i
{m-1\choose i}}(u)
\Biggr)^{{d-1\choose m-1}}
\\
&=&\prod_{m=1}^{d} \prod
_{i=0}^{m-1}g_{d-m+1+i}^{(-1)^i {m-1\choose
i} {d-1\choose m-1}}(u)
\\
&\stackrel{(\ast\ast)} {=}& \prod_{k=1}^d
g_k^{\sum_{m=d+1-k}^d
(-1)^{k+m-d-1} {m-1\choose k+m-d-1} {d-1\choose m-1}}(u)
\\
&=& u.
\end{eqnarray*}
The equality in $(\ast)$ stems from the fact that there are $d-1$ over
$m-1$ shocks $Z^E$ with cardinality $m$ and distribution function
$H_{m,d-m}$ that appear in the stochastic construction of $X_k$. The
equality in $(\ast\ast)$ is yielded by grouping the $g_k$, i.e.
regarding all $g_{d-m+1+i}$ with $d-m+1+i=k$, which is the same as
setting $i=k+m-1-d$. For $k=1$, it is apparent that the exponent in the
second last line is equal to one. For $k\geq2$, due to
\[
\pmatrix{m-1\cr k+m-d-1} \pmatrix{d-1\cr m-1}=\frac{(d-1)!}{(d-k)! (k-1)!} \pmatrix{k-1\cr d-m},
\]
it follows that the exponent equals
\begin{eqnarray*}
&&\frac{(d-1)!}{(d-k)! (k-1)!} \sum_{m=d+1-k}^d
(-1)^{k+m-d-1} \pmatrix{k-1\cr d-m}
\\
&&\quad= \frac{(d-1)!}{(d-k)! (k-1)!} \sum_{m=0}^{k-1}
(-1)^{m} \pmatrix{k-1\cr m}=0.
\end{eqnarray*}
Second, consider the joint distribution function of $(X_1,\ldots,X_d)$.
For $u_1,\ldots,u_d\in(0,1]$, it is given by
%
\begin{equation}\label{joint}
\mathbb{P}(X_k\leq u_k,k=1,\ldots,d)=\prod
_{\varnothing\neq E\subseteq\{
1,\ldots,d\}}\mathbb{P} \bigl(Z^E\leq\min
\{u_k\dvt  k\in E\} \bigr).
\end{equation}
Among all subsets $E$ with cardinality $\vert E\vert=m$, there are
$d-1$ choose $m-1$ elements where $\min\{u_k\dvt k\in E\}=u_{(1)}$.
Analogously, there are $d-k$ choose $m-1$ elements where $\min\{u_l\dvt
l\in E\}=u_{(k)}$, $k\in\{2,\ldots,d-m+1\}$. Thus, \eqref{joint} is
equal to
\begin{eqnarray*}
&&\prod_{m=1}^d \prod
_{k=1}^{d-m+1}\mathbb{P} \bigl(Z^E\leq
u_{(k)}, \vert E\vert=m \bigr)^{{d-k\choose m-1}}
\\
&& \quad= \prod_{m=1}^d \prod
_{k=1}^{d-m+1} H_{m,d-m}^{{d-k \choose
m-1}}(u_{(k)})
\\
&&\quad = \prod_{m=1}^d \prod
_{k=1}^{d-m+1} \prod_{i=0}^{m-1}g_{d-m+1+i}^{(-1)^i {m-1\choose i} {d-k\choose
m-1}}(u_{(k)})
\\
&&\quad= \prod_{k=1}^d \prod
_{m=1}^{d-k+1}\prod_{i=0}^{m-1}g_{d-m+1+i}^{(-1)^i {m-1\choose i} {d-k\choose
m-1}}(u_{(k)})
\\
&&\quad\stackrel{(\ast)} {=} \prod_{k=1}^d
\underbrace{\prod_{n=k}^d
g_n^{\sum_{m=d+1-n}^{d+1-k} (-1)^{m-1+n-d} {m-1\choose m-1+n-d}
{d-k\choose m-1}}(u_{(k)})}_{\mathrm{should} \ \mathrm{be}  \ \mathrm{equal} \ \mathrm{to} \ g_k(u_{(k)})}.
\end{eqnarray*}
Now $(\ast)$ can be derived by sorting all $g_{d-m+1+i}$ with
$d-m+1+i=n$, that is, $i=m-1+n-d$. For $n=k$, it becomes obvious that
the exponent of $g_n$ is equal to one. For $n\in\{k+1,\ldots,d\}$, by
using the same deliberations as for the derivation of the marginal
distributions, the exponent of $g_n$ is given by
\begin{eqnarray*}
&&\frac{(d-k)!}{(d-n)! (n-k)!} \sum_{m=d+1-n}^{d+1-k}(-1)^{m-1+n-d}
\pmatrix{n-k\cr m-1+n-d}
\\
&&\quad= \frac{(d-k)!}{(d-n)! (n-k)!} \sum_{m=0}^{n-k}(-1)^{m}
\pmatrix{n-k\cr m}=0.
\end{eqnarray*}
Summing up, we have
\[
\mathbb{P}(X_k\leq u_k, \forall k=1,\ldots,d)=\prod
_{k=1}^d g_k(u_{(k)})=C(u_1,
\ldots,u_d),
\]
and we can conclude that $C$ is a copula.

$\mbox{(i)}\Rightarrow\mbox{(ii)}$: Let $(U_1,\ldots,U_d)$ be a
random vector with copula $C$ in \eqref{G} as distribution function.
Moreover, assume for a moment that $g_k, k=2,\ldots,d$, is strictly
positive on $(0,1]$. Then, for $u,v\in(0,1], u<v$, $G_{j,k}$ has the
representation
\begin{eqnarray*}
G_{j,k}(u,v)&=&\frac{1}{\prod_{m=1}^{k}g_m(u)} \Biggl(\mathbb {P}(A_{\varnothing})-
\sum_{i=1}^{j}(-1)^{i+1}\mathop{
\sum_{L\subseteq
\{k+1,\ldots,k+j\}\dvt}}_{\vert L\vert=i}\mathbb{P} \biggl(\bigcap
_{l\in L} A_l \biggr) \Biggr),
\\
A_l &:=& \biggl(\bigcap_{m\in\{1,\ldots,k,l\}}
\{U_m\leq u\} \biggr)\cap \biggl(\bigcap
_{m\in\{k+1,\ldots,k+j\}\setminus\{l\}} \{U_m\leq v\} \biggr),
\\
A_{\varnothing}&:=& \{U_1\leq u,\ldots,U_k\leq
u,U_{k+1}\leq v,\ldots ,U_{k+j}\leq v\}.
\end{eqnarray*}
Applying the
principle of inclusion and exclusion (see \cite{Billingsley1995}, page 24), we have
\begin{eqnarray*}
G_{j,k}(u,v)&=& \frac{1}{\prod_{m=1}^{k}g_m(u)} \Biggl(\mathbb {P}(A_{\varnothing})-
\mathbb{P} \Biggl(\bigcup_{l=k+1}^{k+j}
A_l \Biggr) \Biggr)
\\[-2pt]
&=& \frac{1}{\prod_{m=1}^{k}g_m(u)} \mathbb{P}(A)\geq0,
\end{eqnarray*}
where\vspace*{-3pt}
\[
A:= \bigl\{U_1\leq u,\ldots,U_k\leq
u,U_{k+1}\in[u,v],\ldots,U_{k+j}\in [u,v] \bigr\}.
\]
Strict positivity of $g_k$ (which we have assumed so far) as well as
continuity on $(0,1]$ for $k=2,\ldots,d$ can be shown by induction. To
begin with, assume that there is a $u^\ast:=\sup\{u\geq0 \dvt g_2(u)=0\}>0$.
As $C$ is a copula and hence continuous, it follows that $g_2(u^\ast
)=0$, such that for $v>u^\ast$,
\[
\mathbb{P} \bigl(U_1\in\bigl[u^{\ast},v
\bigr],U_2\in\bigl[u^{\ast},v\bigr] \bigr)=G_{2,0}
\bigl(u^\ast ,v\bigr)=v g_2(v)-2 u^\ast
g_2(v)<0
\]
for $v$ sufficiently close to $u^\ast$. This is a contradiction and
hence $g_2(u)>0$ for $u\in(0,1]$. Similarly, to show continuity, assume
that  there is a $v^\ast\in(0,1]$ such that $g_2(v^\ast-):=\lim_{u\nearrow v^\ast} g_2(u)<g_2(v^\ast)$. Then
\begin{eqnarray*}
0 &\leq & \lim_{u\nearrow v^\ast} G_{2,0}\bigl(u,v^\ast
\bigr)=\lim_{u\nearrow v^\ast} \bigl(v^\ast g_2
\bigl(v^\ast\bigr)-2 u g_2\bigl(v^\ast\bigr)+u
g_2(u) \bigr)
\\
&=&-v^\ast g_2\bigl(v^\ast\bigr)+v^\ast
g_2\bigl(v^\ast-\bigr)<0,
\end{eqnarray*}
which is a contradiction. Hence, there is no such $v^\ast$ and $g_2$ is
left-continuous on $(0,1]$. Analogously, if $g_2(u^\ast+):=\lim_{v\searrow u^\ast}>g_2(u^\ast)$ for an $u^\ast\in(0,1)$,
$G_{2,0}(u^\ast,v)$ becomes negative for sufficiently small $v>u^\ast$.
Consequently, $g_2$ is continuous on $(0,1]$.

For the induction step $k-1\mapsto k$, note that
\[
G_{2,k-1}(u,v)=\frac{1}{\prod_{m=1}^{k-1}g_m(u)} \mathbb{P} \bigl(U_1\leq
u,\ldots,U_{k-1}\leq u,U_k\in[u,v],U_{k+1}
\in[u,v] \bigr)
\]
induces $0\leq G_{2,k-1}(u,v)=g_{k-1}(v) g_k(v)-2 g_{k-1}(u)
g_k(v)+g_{k-1}(u) g_k(u)$. By the same arguments as for the induction
start, this implies that $g_k$ is both strictly positive and continuous
on $(0,1]$.

\begin{rem}[(Alternative interpretation of $G_{j,k}$)]\label{XY}
In a very similar way, one can show that if $H_{m,k+j-m}$ in Theorem~\ref{maintheorem}(iii) is an element of $\mathcal{D}$
for all $m\in\{1,\ldots,j\}$, $G_{j,k}$ can be expressed as\vspace*{-6pt}
\[
G_{j,k}(u,v) = \mathbb{P}\bigl(X_{k+1}\in[u,v],
\ldots,X_{k+j}\in[u,v]\bigr),
\]
where\vspace*{-3pt}
\[
X_{l}:=\max\bigl\{Z^E\dvt  k\in E\bigr\}, l=k+1,\ldots,k+j  \mbox{ and }
Z^E\sim H_{m,k+j-m}   \mbox{ for } \vert E\vert=m,
\]
with independent random variables $Z^E, \varnothing\neq E\subseteq\{
k+1,\ldots,k+j\}$. However, note that in the proof of ``$\mbox{(i)}\Rightarrow
\mbox{(ii)}$'' above, we solely require that $C$ is a copula and do not assume
increasingness of $H_{m,k+j-m}$, which is why we cannot apply the
alternative interpretation in the present case.
\end{rem}

$\mbox{(ii)}\Rightarrow\mbox{(iii)}$  Let $G_{j,k}(u,v)\geq0$ for
all $0< u< v\leq1, k\in\mathbb{N}_{0}, j\in\mathbb{N}$ with $k+j\leq
d$. By the proof of ``$\mbox{(i)}\Rightarrow\mbox{(ii)}$'' above, this implies $g_k$,
$k=2,\ldots,d$ to be strictly positive and continuous on $(0,1]$. The
idea of this part of the proof is to establish a connection between
$G_{j,k}$ and a related stochastic model similar to Remark~\ref{XY} in
order to derive reasonable estimates that help to derive the required
conditions in $\mbox{(iii)}$. For readability, we are going to proceed by
induction.

Suppose that we have already shown that for a $j-1\in\{1,\ldots,d\}$,
the conditions $G_{i,k}(u,v)\geq0$ for all $1\leq i\leq j-1$ and
$k+i\leq d$ imply that $H_{i,k}$ is increasing for all $1\leq i\leq
j-1$ and $k+i\leq d$. For $j=2$, this is obviously satisfied and the
induction basis is established. In order to carry out the induction
step, we need to show that $H_{j,k}$ is increasing for all $k+j\leq d$
($H_{i,k}$ for $i\leq j-1$ are increasing by induction hypothesis).
This is shown in several steps.
\begin{longlist}
\item[Step 1 (Main observation):] There is a useful decomposition of $G_{j,k}$
that we require below.
\end{longlist}

\begin{lem}[(Decomposition of $G_{j,k}$)]\label{decomposition}
Instead of $G_{j,k}$, write $G_{g_{k+1},\ldots,g_{k+j}}$ to emphasize
the dependence of $G_{j,k}$ on the functions $g_{k+1},\ldots,g_{k+j}$.
It holds that
%
\begin{eqnarray}
G_{g_{k+1},\ldots,g_{k+j}}(u,v)&=& \tilde{g}_{k+1}(v)
g_{k+2}(v) \ldots g_{k+j}(v) \biggl(\frac{g_{k+1}(v)}{\tilde{g}_{k+1}(v)}-
\frac
{g_{k+1}(u)}{\tilde{g}_{k+1}(u)} \biggr)
\nonumber
\\[-9pt]
\label{G2}
\\[-9pt]
\nonumber
&&{}+\frac{g_{k+1}(u)}{\tilde{g}_{k+1}(u)} G_{\tilde
{g}_{k+1},g_{k+2},\ldots,g_{k+j}}(u,v),\qquad 0<u<v\leq1,
\end{eqnarray}
%
for an arbitrary function $\tilde{g}_{k+1}$ that is unequal to zero on $(0,1)$.
\end{lem}

\begin{pf}
The decomposition consists of nothing else than changing the last
summand of $G_{g_{k+1},\ldots,g_{k+j}}$ (ending up with the last line
in equation \eqref{G2}) and adding the resulting difference as an extra
term (corresponding to the first line in equation \eqref{G2}). The
non-zero condition for $\tilde{g}_{k+1}$ is required for well-defined quotients.
%
\end{pf}

Define $\tilde{g}_{k+1}:=g_{k+1}/H_{j,k}$, which is
continuous and strictly positive on $(0,1]$ as seen earlier, and note
that Lemma~\ref{decomposition} then yields
%
\begin{eqnarray}
0 &\leq &  G_{g_{k+1},\ldots,g_{k+j}}(u,v)
\nonumber
\\[-9pt]
\label{decomp}
\\[-9pt]
\nonumber
&=&\tilde{g}_{k+1}(v) g_{k+2}(v)
\cdots g_{k+j}(v) \bigl(H_{j,k}(v)-H_{j,k}(u) \bigr)
+H_{j,k}(u) G_{\tilde{g}_{k+1},g_{k+2},\ldots,g_{k+j}}(u,v).
\end{eqnarray}
We want to conclude that $H_{j,k}(v)\geq H_{j,k}(u)$. Therefore, we
have to prove that the second summand is not responsible for
non-negativity of $G_{g_{k+1},\ldots,g_{k+j}}$. The crucial consequence
of \eqref{decomp} is that $G_{\tilde{g}_{k+1},g_{k+2},\ldots,g_{k+j}}$
can be related to a stochastic model. To this end, we want to apply
Remark~\ref{XY} to $G_{\tilde{g}_{k+1},g_{k+2},\ldots,g_{k+j}}$. In
order to do so, one has to make sure that the corresponding functions
$\tilde{H}_{m,k+j-m}$ (which are defined just like $H_{m,k+j-m}$,
however with replacing $g_{k+1}$ by $\tilde{g}_{k+1}$) are distribution
functions in $\mathcal{D}$ for all $m=1,\ldots,j$. Due to the
definition of $\tilde{g}_{k+1}$, it holds that\vspace*{-2pt}
\[
\tilde{H}_{m,k+j-m}=
\cases{ H_{m,k+j-m}, &\quad\mbox{for
$m=1,\ldots,j-1$},\vspace*{2pt}
\cr
1, &\quad\mbox{for $m=j$}.}
\]
As $H_{m,k+j-m}$ are distribution functions for $m=1,\ldots,j-1$ by
induction hypothesis and $\tilde{H}_{m,k+j-m}$ is a degenerated
distribution function for $m=j$, the requirements of Remark~\ref{XY}
are satisfied.
\begin{longlist}
\item[Step 2 (Stochastic model):] As a consequence of Remark~\ref{XY},
\end{longlist}
\[
G_{\tilde{g}_{k+1},g_{k+2},\ldots,g_{k+j}}(u,v)=\mathbb{P} \bigl(X_{k+1},\ldots,X_{k+j}
\in[u,v] \bigr),
\]
where\vspace*{-3pt}
\begin{eqnarray*}
X_l &:=& \max \bigl\{Z^E, E\subset\{k+1,\ldots,k+j\}, E
\cap\{l\}\neq\varnothing \bigr\}, l=k+1,\ldots,k+j,
\\
Z_E &\sim & H_{m,k+j-m} \mbox{ for $\vert E\vert=m$},
\end{eqnarray*}
with independent random variables $Z^E, \varnothing\neq E\subset\{
k+1,\ldots,k+j\}$. Thus,
\begin{eqnarray*}
&&\mathbb{P} \bigl(X_{k+1},\ldots,X_{k+j}\in[u,v] \bigr)
\\[-2pt]
&&\quad=\mathbb{P} \Biggl(\underbrace{\bigcap_{l=k+1}^{k+j}
\bigl\{\max \bigl\{Z^E, E\cap\{l\}\neq\varnothing \bigr\}\in[u,v] \bigr
\}}_{:=A} \Biggr)
\end{eqnarray*}
requires that all $Z^E$ are less than or equal to $v$ and -- as there
is no common shock with $\vert E\vert=j$ due to $\tilde
{H}_{j,k}(x)\equiv1$ -- at least two $Z^I, Z^J, I,J\subset\{k+1,\ldots
,k+j\},I\neq J$, need to be in the interval $[u,v]$. This implies that
\begin{eqnarray*}
A\subset\mathop{\bigcup_{\varnothing\neq I,J\subset\{k+1,\ldots,k+j\}}}_{I\neq J}
\bigl\{u \leq Z^I,Z^J\leq v \bigr\}.
\end{eqnarray*}
Moreover,\vspace*{1pt} as $\mathbb{P}(\bigcup_{i=1}^n A_i)\leq\sum_{i=1}^n \mathbb
{P}(A_i)$ for arbitrary $A_i\in\mathcal{F}$, and as there are $(2^j-2)$
choose~$2$ possibilities to pick $Z^I,Z^J\in[u,v]$, we have
%
\begin{eqnarray}
G_{\tilde{g}_{k+1},g_{k+2},\ldots,g_{k+j}}(u,v)&=&\mathbb{P}(A)\leq\mathop{\sum
_{\varnothing\neq I,J\subset\{k+1,\ldots,k+j\}}}_{I\neq
J}\mathbb{P} \bigl( \bigl\{u \leq
Z^I,Z^J\leq v \bigr\} \bigr)
\nonumber
\\[-9pt]
\label{sum}
\\[-9pt]
\nonumber
&\leq & \underbrace{\pmatrix{2^j-2\cr 2}}_{:=\,b} \max
_{m=1,\ldots
,j-1} \bigl\{ \bigl(H_{m,k+j-m}(v)-H_{m,k+j-m}(u)
\bigr)^2 \bigr\}.
\end{eqnarray}
%
%

\begin{longlist}
\item[Step 3 (Lipschitz-continuity):] Using equation \eqref{sum}, we are going
to derive Lipschitz-continuity-type results for $G_{\tilde
{g}_{k+1},g_{k+2},\ldots,g_{k+j}}$. In order to do so, the following
lemma is helpful.
\end{longlist}

\begin{lem}\label{Lipschitz} For $k\in\mathbb{N}_0, j\geq2$, let
$H_{1,k},\ldots,H_{j,k}\dvtx (0,1]\rightarrow(0,1]$ and $H_{1,k+1},\ldots
,\break H_{j-1,k+1}\dvtx [0,1]\rightarrow[0,1]$ be increasing functions with
$H_{l,k}=H_{l-1,k}/H_{l-1,k+1}$ for $l\in\{2,\ldots,j\}$. Then it holds that
\[
0\leq H_{j,k}(v)-H_{j,k}(u)\leq \Biggl(\prod
_{l=1}^{j-1}\frac
{1}{H_{l,k+1}(u)} \Biggr)
\bigl(H_{1,k}(v)-H_{1,k}(u) \bigr).
\]
\begin{pf}
For $j=2$ and $k\in\mathbb{N}_0$, we have
\begin{eqnarray*}
0& \leq&H_{2,k}(v)-H_{2,k}(u)=\frac{1}{H_{1,k+1}(u)} \biggl(
\underbrace {\frac{H_{1,k+1}(u)}{H_{1,k+1}(v)}}_{\leq1} H_{1,k}(v)-H_{1,k}(u)
\biggr)
\\
&\leq&\frac{1}{H_{1,k+1}(u)} \bigl(H_{1,k}(v)-H_{1,k}(u) \bigr).
\end{eqnarray*}
For $j\mapsto j+1$, the claim follows by simple induction.
\end{pf}
\end{lem}

Applying Lemma~\ref{Lipschitz} to equation \eqref{sum},
$G_{\tilde{g}_{k+1},\ldots,g_{k+j}}(u,v)$ has an upper bound
%
\begin{equation}\label{sumg}
b \cdot \max_{m=1,\ldots,j-1} \Biggl\{ \Biggl(\prod
_{l=1}^{m-1}\frac
{1}{H_{l,k+j-m+1}(u)} \Biggr)^2
\bigl(\underbrace {g_{k+j-m+1}(v)-g_{k+j-m+1}(u)}_{=\,H_{k+j-m}(v)\,-\,H_{k+j-m}(u)}
\bigr)^2 \Biggr\}.
\end{equation}
This can be further simplified as
\begin{eqnarray*}
0&\leq & G_{2,0}(u,v)=g_1(v) g_2(v)-2
g_1(u) g_2(v)+g_1(u) g_2(u)
\\
&=& g_2(v) \bigl(g_1(v)-g_1(u)
\bigr)-g_1(u) \bigl(g_2(v)-g_2(u) \bigr)
\\
&& \hspace*{-5pt}\Longleftrightarrow\quad  g_2(v)-g_2(u)\leq  \frac{g_2(v)}{g_1(u)}
\bigl(g_1(v)-g_1(u) \bigr)=\frac{g_2(v)}{u} (v-u )
\end{eqnarray*}
and since $g_1$ is the identity by definition. Analogously, one has
\[
0\leq G_{2,k}(u,v) \quad \Longleftrightarrow \quad g_{k+2}(v)-g_{k+2}(u)
\leq\frac
{g_{k+2}(v)}{g_{k+1}(u)} \bigl(g_{k+1}(v)-g_{k+1}(u) \bigr).
\]
By induction over $k$, one can conclude that
\[
g_k(v)-g_k(u)\leq \Biggl(\prod
_{l=1}^{k-1}\frac{g_{l+1}(v)}{g_l(u)} \Biggr) (v-u)
\qquad \mbox{for all $k$ in concern}.
\]
Applying this result, the expression in \eqref{sumg} is less than or
equal to
\begin{eqnarray*}
&& b \cdot \max_{m=1,\ldots,j-1} \Biggl\{ \Biggl(\prod
_{l=1}^{m-1}\frac
{1}{H_{l,k+j-m+1}(u)} \Biggr)^2
\Biggl(\prod_{l=1}^{k+j-m}\frac
{g_{l+1}(v)}{g_l(u)}
\Biggr)^2 (v-u)^2 \Biggr\}
\\
&&\quad = p_{j,k}(u,v) (v-u)^2,
\end{eqnarray*}
with
\[
p_{j,k}(u,v):=b \cdot \max_{m=1,\ldots,j-1} \Biggl\{ \Biggl(
\prod_{l=1}^{m-1}\frac{1}{H_{l,k+j-m+1}(u)}
\Biggr)^2 \Biggl(\prod_{l=1}^{k+j-m}
\frac{g_{l+1}(v)}{g_l(u)} \Biggr)^2 \Biggr\}.
\]
Additionally, due to the monotonicity of the $g_k$ and $H_{j,k}$
appearing in $p_{j,k}$, we can conclude that for any $u_0,v_0\in(0,1],
u_0<v_0$, it holds that $p_{j,k}(u,v)\leq p_{j,k}(u_0,v_0)$ for all
$u,v\in[u_0,v_0]$, \mbox{$u\leq v$}. Combining all those observations, one ends
up with
%
\begin{equation}\label{finallipschitz}
0\leq G_{\tilde{g}_{k+1},g_{k+2},\ldots,g_{k+j}}(u,v)\leq p_{j,k}(u_0,v_0)
(v-u)^2, \qquad u,v\in[u_0,v_0], u\leq v.
\end{equation}

Step 4 (Proof by contradiction): Finally, we can proceed similarly to
the proof in the bivariate case depicted in \cite{Durante2008}, page 67.
Assume that $H_{j,k}$ is not increasing and that there exist $u_0,v_0\in
(0,1], u_0<v_0$, such that
%
\[
H_{j,k}(v_0)-H_{j,k}(u_0)=-a(u_0,v_0)
(v_0-u_0), \qquad a(u_0,v_0)>0.
\]
Consequently, by continuity of the $g_k$ and hence
$H_{j,k}$, for every $\varepsilon\in(0,v_0-u_0]$,
there are $u_\varepsilon, v_\varepsilon\in[u_0,v_0],
u_\varepsilon=v_\varepsilon-\varepsilon$, such that
\begin{equation}\label{k}
H_{j,k}(v_\varepsilon)-H_{j,k}(u_\varepsilon)
\leq-a(u_0,v_0) (v_\varepsilon -u_\varepsilon)=-a(u_0,v_0)
\varepsilon.
\end{equation}
Independently of this assumption, we can split the positive and
negative powers in $H_{j,k}$, yielding
\[
H_{j,k}(u)=\prod_{i=0}^{j-1}g_{k+1+i}^{(-1)^i {j-1 \choose
i}}(u)=
\frac{\prod_{i=0}^{\lfloor ({j-1})/{2} \rfloor}g_{k+1+2
i}^{{j-1 \choose2 i}}(u)}{\prod_{i=0}^{\lfloor ({j-2})/{2} \rfloor
}g_{k+1+2 i+1}^{{j-1 \choose2 i+1}}(u)},\qquad u>0,
\]
with ``$\lfloor\cdot\rfloor$'' denoting the floor function, such that
for $u\in[u_0,v_0]$, it holds by the monotonicity of the $g_k$ that
\begin{eqnarray*}
p_{\min}(u_0,v_0)&:=& \frac{\prod_{i=0}^{\lfloor ({j-1})/{2} \rfloor
}g_{k+1+2 i}^{{j-1 \choose2 i}}(u_0)}{\prod_{i=0}^{\lfloor ({j-2})/{2} \rfloor}g_{k+1+2 i+1}^{{j-1 \choose2 i+1}}(v_0)}\leq
H_{j,k}(u)
\\
&\leq& \frac{\prod_{i=0}^{\lfloor ({j-1})/{2} \rfloor}g_{k+1+2
i}^{{j-1 \choose2 i}}(v_0)}{\prod_{i=0}^{\lfloor ({j-2})/{2}
\rfloor}g_{k+1+2 i+1}^{{j-1 \choose2 i+1}}(u_0)}=:p_{\max}(u_0,v_0).
\end{eqnarray*}
Plugging $u_\varepsilon, v_\varepsilon$ into equation \eqref{decomp} and
using all previous results yields
\begin{eqnarray*}
0 &\leq &  G_{g_{k+1},\ldots,g_{k+j}}(u_\varepsilon,v_\varepsilon)=\underbrace {
\tilde{g}_{k+1}(v_\varepsilon)}_{=\frac{g_{k+1}(v_\varepsilon
)}{H_{j,k}(v_\varepsilon)}} g_{k+2}(v_\varepsilon)
\cdots g_{k+j}(v_\varepsilon) \bigl(\underbrace{H_{j,k}(v_\varepsilon
)-H_{j,k}(u_\varepsilon)}_{\leq-a(u_0,v_0) \varepsilon} \bigr)
\\
&&{}+H_{j,k}(u_\varepsilon) G_{\tilde{g}_{k+1},\ldots,g_{k+j}}(u_\varepsilon
,v_\varepsilon)
\\
&\stackrel{\mathrm{\eqref{finallipschitz}}} {\leq} &\frac{g_{k+1}(u_0)}{p_{\max
}(u_0,v_0)}
g_{k+2}(u_0) \cdots g_{k+j}(u_0)
\bigl(-a(u_0,v_0) \varepsilon \bigr)
\\
&&{}+p_{\max}(u_0,v_0) p_{j,k}(u_0,v_0)
\varepsilon^2.
\end{eqnarray*}
Thus, for sufficiently small $\varepsilon$, $G_{j,k}$ becomes negative and
yields a contradiction. Consequently, $H_{j,k}$ has to be increasing
and the induction is complete.

\begin{rem}The proofs ``$\mbox{(i)}\Rightarrow\mbox{(ii)}$'' and ``$\mbox{(ii)}\Rightarrow
\mbox{(iii)}$'' generalize parts of the proof ideas for Theorem~1.1 in \cite
{Mai2014}. The major generalization consists in omitting
differentiability of the $g_k$, which substantially complicates the
calculations and requires the alternative proof techniques picked up in
``$\mbox{(ii)}\Rightarrow\mbox{(iii)}$'' above.
\end{rem}

\begin{longlist}
\item[$\mbox{(iii)}\Rightarrow\mbox{(iv)}$:] Trivial, as $\mbox{(iv)}$ is a special case
of $\mbox{(iii)}$ for $j=m, k=d-m$.\hfill\qed
\end{longlist}
\noqed\end{pf*}

\section{}\label{appsubclass}
\begin{pf*}{Proof of Proposition \protect\ref{subclass}}
The survival function of each $X_k$ is given by
\begin{eqnarray*}
\bar{F}_1(x)&:=& \mathbb{P}(X_k>x)=\mathbb{P}(E_k>
\Lambda_x)=\mathbb {E} \bigl[\mathbb{P}(E_k>
\Lambda_x\vert\Lambda_x) \bigr]
\\
&=&\mathbb{E} \bigl[\mathrm{e}^{-\Lambda_x} \bigr]=\mathrm{e}^{-\Psi_x(1)},\qquad x\geq0.
\end{eqnarray*}
The joint survival function of $(X_1,\ldots,X_d)$ can be derived analogously.
For $\textbf{x}:=(x_1,\ldots, x_d)\geq0$, with the convention $x_0:=0$,
it is given by
\begin{eqnarray*}
\bar{F}_d(\mathbf{x})&:=& \mathbb{P}(X_1>x_1,
\ldots,X_d>x_d)=\mathbb {E} \bigl[\mathrm{e}^{-\sum_{k=1}^d \Lambda_{x_k}} \bigr]
\\
&=& \mathbb{E} \bigl[\mathrm{e}^{-\sum_{k=1}^d (d-k+1) (\Lambda_{x_{(k)}}-\Lambda
_{x_{(k-1)}})} \bigr]
=\prod_{k=1}^d \mathbb{E}
\bigl[\mathrm{e}^{- (d-k+1) (\Lambda_{x_{(k)}}-\Lambda
_{x_{(k-1)}})} \bigr]
\\
&=&\prod_{k=1}^d \exp \bigl(-
\Psi_{x_{(k)}}(d-k+1)+\Psi _{x_{(k-1)}}(d-k+1) \bigr)
\\
&=&\prod_{k=1}^d \exp \bigl(-
\Psi_{x_{(k)}}(d-k+1)+\Psi_{x_{(k)}}(d-k) \bigr). 
\end{eqnarray*}
%
Due to the stochastic continuity of $\Lambda$, $x\mapsto\Psi_x(1)$ is
continuous. Thus, the unique survival copula $C_{\Lambda,d}$ of
$(X_1,\ldots,X_d)$ is defined as
%
\begin{eqnarray}
C_{\Lambda,d}(u_1,\ldots,u_d)&:=&
\bar{F}_d \bigl(\bar{F}_1^{-1}(u_1),
\ldots ,\bar{F}_1^{-1}(u_d) \bigr)
\nonumber
\\
\label{Cadditive}
&=&\prod_{k=1}^d \exp \bigl(-
\Psi_{\bar{F}_1^{-1}(u_{(k)})}(k)+\Psi_{\bar
{F}_1^{-1}(u_{(k)})}(k-1) \bigr)
\\
&=&\prod_{k=1}^d g_k(u_{(k)}),\nonumber
\end{eqnarray}
where
\[
g_k(u):=\exp \bigl(-\Psi_{\bar{F}_1^{-1}(u)}(k)+\Psi_{\bar
{F}_1^{-1}(u)}(k-1)
\bigr).
\]
By construction, $g_1=\mbox{id}_{[0,1]}$ and $g_2(1)=\cdots=g_d(1)=1$.
\end{pf*}


%
%
\end{appendix}

\section*{Acknowledgements}
We are grateful to the anonymous referees for their valuable
suggestions and remarks on a previous version of this manuscript.






\printhistory
\end{document}